\documentclass[11pt]{article}
\usepackage[margin=1in]{geometry}
\usepackage[T1]{fontenc}
\usepackage{lmodern}
\usepackage{amsmath,amssymb,amsthm,mathtools}
\usepackage{microtype}
\usepackage[colorlinks=true,linkcolor=blue,citecolor=blue,urlcolor=blue]{hyperref}
\hypersetup{
  pdftitle={The Order of the Inertia Bound for Random Graphs},
  pdfauthor={Zhengbo Chen and Xiao-Dong Zhang}
}
\numberwithin{equation}{section}
\newtheorem{theorem}{Theorem}[section]
\newtheorem{conjecture}[theorem]{Conjecture}
\newtheorem{corollary}[theorem]{Corollary}

\newtheorem{lemma}[theorem]{Lemma}
\theoremstyle{definition}
\newtheorem{definition}[theorem]{Definition}
\newtheorem{question}[theorem]{Question}
\theoremstyle{remark}
\newtheorem{remark}[theorem]{Remark}
\DeclareMathOperator{\tr}{tr}
\DeclareMathOperator{\rank}{rank}
\DeclareMathOperator{\diag}{diag}
\DeclareMathOperator{\Diag}{Diag}
\newcommand{\one}{\mathbf 1}
\newcommand{\C}{\mathbb C}
\newcommand{\R}{\mathbb R}
\newcommand{\W}{\mathcal W_{\mathrm H}}
\title{The Order of the Inertia Bound for Random Graphs}
\author{Zhengbo Chen \quad Xiao-Dong Zhang\thanks{Corresponding author.
Email: \href{mailto:xiaodong@sjtu.edu.cn}{\nolinkurl{xiaodong@sjtu.edu.cn}}.}\\[0.5em]
\small School of Mathematical Sciences, MOE-LSC, SHL-MAC\\
\small Shanghai Jiao Tong University, Shanghai 200240, China}
\date{20 November 2026}
\begin{document}
\maketitle

\begin{abstract}
For an $n$-vertex graph $G$, let $n_{\geq0}(G)$ be the minimum number
of nonnegative eigenvalues among zero-diagonal Hermitian matrices
supported on $G$, allowing zero edge weights. For $n\geq2$ and
$G\sim G(n,1/2)$, we prove $n_{\geq0}(G)\geq c_2n/\log(2n)$
with probability at least $1-\exp(-c_1n\log(2n))$, where
$c_1,c_2>0$ are absolute constants. This proves a Hermitian extension
of the Kwan--Wigderson conjecture. With high probability,
$n_{\geq0}(G)=\Theta(n/\log n)$ and
$n_{\geq0}(G)/\alpha(G)=\Theta(n/(\log n)^2)$, where $\alpha(G)$
is the independence number.

Using Alon's constructions, we obtain, for every sufficiently large
$n$, an $n$-vertex graph $G$ with $\alpha(G)=2$ and
$n_{\geq0}(G)=\Omega(n^{1/3})$. For infinitely many $n$, we also
obtain $n$-vertex graphs $H$ with Shannon capacity $3$ and
$n_{\geq0}(H)=\Omega(n^{1/4})$.
\end{abstract}

\section{Introduction}\label{sec:introduction}

Let $G$ be a finite simple graph on $n$.
Write $A_G$ for its adjacency matrix, $\overline G$ for its complement,
$\alpha(G)$ for its independence number, and $\chi(G)$ for its
chromatic number. The independence number is the maximum size of a
set of pairwise nonadjacent vertices. Computing $\alpha(G)$ is
NP-hard, motivating spectral upper bounds such as the inertia bound.

Let $\W(G)$ be the set of Hermitian matrices
$A=(a_{ij})\in\C^{n\times n}$ with zero diagonal and $a_{ij}=0$
whenever $ij\notin E(G)$. These are the Hermitian weighted adjacency
matrices of $G$, with zero weights permitted on edges.
Throughout, $\log$ denotes the natural logarithm, and an event
holds with high probability if its probability tends to one as
$n\to\infty$.

For a Hermitian matrix $A$, let $n_+(A)$, $n_-(A)$, and $n_0(A)$
count its positive, negative, and zero eigenvalues, respectively,
with multiplicity. Define its nonnegative inertia index by
$n_{\geq0}(A)\triangleq n_+(A)+n_0(A)$, and the optimal inertia
bound of $G$ by
\begin{equation}\label{eq:inertia-definition}
 n_{\geq0}(G)\triangleq\min_{A\in\W(G)}n_{\geq0}(A).
\end{equation}
Since $A$ and $-A$ are both admissible, this is also the minimum
of the usual weighted inertia bound $n_0(A)+\min\{n_+(A),n_-(A)\}$.
For a prescribed Hermitian matrix with rational real and imaginary
parts, inertia is computable in polynomial time by exact linear
algebra, whereas optimization over the weights is a separate
problem~\cite[Section~1]{KW}.

The following classical inertia bound is recalled
in~\cite[Theorem~1.2]{KW}; see~\cite[equation~(1.1)]{TZE}
for the Hermitian formulation.
\begin{theorem}[Inertia bound]\label{thm:inertia-bound}
Let $G$ be a finite simple graph. For every $A\in\W(G)$,
\[
 \alpha(G)\leq n_0(A)+\min\{n_+(A),n_-(A)\}
 \leq n_{\geq0}(A).
\]
Consequently,
\begin{equation}\label{eq:alpha-inertia}
 \alpha(G)\leq n_{\geq0}(G).
\end{equation}
\end{theorem}

To measure how far this bound can be from equality, we consider
both the additive and multiplicative gaps.
\begin{question}\label{ques:inertia-gap}
For graphs $G$ on $n\geq1$ vertices, how large can the difference
and the ratio
\[
 n_{\geq0}(G)-\alpha(G)
 \qquad\text{and}\qquad
 \frac{n_{\geq0}(G)}{\alpha(G)}
\]
be? How do these quantities behave for typical graphs or for graphs
with bounded independence number?
\end{question}

A basic question is whether suitable weights always make the inertia
bound equal to the independence number. Godsil posed this question
for real symmetric weights; see~\cite[Question~1.3]{KW}.
Sinkovic~\cite{Sinkovic} answered it negatively: the Paley graph on
$17$ vertices has independence number $3$, whereas every real
symmetric weighted adjacency matrix of this graph has at least $4$
nonnegative eigenvalues.

Kwan and Wigderson~\cite[Theorem~1.6]{KW} subsequently constructed,
for infinitely many $n$, regular graphs satisfying
$\alpha(G)\leq4n^{3/4}$ and $n_{\geq0}(G)\geq n/4$.
Consequently, these graphs have an additive gap of $\Omega(n)$
and a multiplicative gap of $\Omega(n^{1/4})$.
Tang, Zhang, and Elphick~\cite[Theorem~1.2]{TZE} proved that
$\alpha(G)\leq2$ can coexist with $n_{\geq0}(G)=\Omega(n^{1/9})$,
thereby resolving Kwan--Wigderson's Conjecture~4.2.
They also obtained multiplicative gaps of size $n^{1-o(1)}$
for other graph families~\cite[Proposition~1.5]{TZE}.

These results show that the gap can be large even when the
independence number is bounded. It is therefore also natural to
ask which graph classes admit equality.
\begin{question}\label{ques:inertia-equality}
For which classes of graphs does the equality
\[
 n_{\geq0}(G)=\alpha(G)
\]
hold for every graph $G$ in the class?
\end{question}

Perfect graphs form one such class; see also~\cite[Section~1]{Sinkovic}.
To see this, partition $V(G)$ into cliques and take the adjacency
matrix of their disjoint union. This matrix belongs to $\W(G)$.
Each clique contributes exactly one nonnegative eigenvalue, including
the case of a singleton. A partition with $\chi(\overline G)$ parts thus gives
\begin{equation}\label{eq:inertia-clique-cover}
 \alpha(G)\leq n_{\geq0}(G)\leq\chi(\overline G).
\end{equation}
For perfect graphs, $\chi(\overline G)=\alpha(G)$, so this inequality
gives $n_{\geq0}(G)=\alpha(G)$.

We address Question~\ref{ques:inertia-gap} for random graphs and
for graphs with independence number two. We develop a method for
proving inertia lower bounds that hold for every admissible weight
matrix. It determines the typical order of
$n_{\geq0}(G)/\alpha(G)$ for $G\sim G(n,1/2)$ and gives graphs with
$\alpha(G)=2$ and $n_{\geq0}(G)=\Omega(n^{1/3})$.

\subsection{The gap between \texorpdfstring{$n_{\geq0}(G)$}{n>=0(G)} and \texorpdfstring{$\alpha(G)$}{alpha(G)} in random graphs}

For a typical graph, can a suitable choice of weights make the
inertia bound close to the independence number? Ihringer asked
whether equality holds for almost all graphs, as discussed
in~\cite[Section~4]{KW}. Kwan and
Wigderson~\cite[Conjecture~4.1]{KW} proposed the following conjecture,
which predicts a substantial separation for $G\sim G(n,1/2)$, whose
independence number is only of order $\log n$.
\begin{conjecture}
\label{conj:KW}
Let $G\sim G(n,1/2)$. With probability $1-o(1)$ as $n\to\infty$,
every real symmetric weighted adjacency matrix $A$ of $G$ satisfies
\begin{equation}\label{eq:KW-conjecture}
 n_{\geq0}(A)=\Omega(n/\log n),
\end{equation}
where the implied constant is absolute.
\end{conjecture}

We establish this bound for Hermitian weights, thereby proving the
conjecture with an explicit failure probability. Together with the
known estimates for clique partitions and independence numbers,
this quantifies the unavoidable loss in estimating the independence
number by inertia.
\begin{theorem}\label{thm:main-random}
Let $G\sim G(n,1/2)$.
\begin{enumerate}
\item\label{item:random-lower}
There are absolute constants $c_1,c_2>0$ such that, for every $n\geq2$,
with probability at least $1-\exp(-c_1n\log(2n))$, every Hermitian
weighted adjacency matrix $A$ of $G$ satisfies
\[
 n_{\geq0}(A)\geq c_2\frac{n}{\log(2n)}.
\]
\item\label{item:random-gap}
With high probability,
\begin{equation}\label{eq:random-inertia-gap}
 \frac{n_{\geq0}(G)}{\alpha(G)}
 =\Theta\!\left(\frac{n}{(\log n)^2}\right).
\end{equation}
\end{enumerate}
\end{theorem}
We prove part~\ref{item:random-lower} in Section~\ref{sec:random-energy}
and derive part~\ref{item:random-gap} below.
The event in Theorem~\ref{thm:main-random}(\ref{item:random-lower})
depends only on $G$ and works for every admissible choice of weights, including choices
made after the graph is revealed. Thus even the optimal inertia
bound is at least of order $n/\log n$; the separation from
$\alpha(G)$ cannot be removed by a better choice of weights.

To determine the size of this separation, we also need an upper
bound on the optimal inertia bound. The clique-partition observation
of Alon, recorded in~\cite[Section~4]{KW}, gives precisely
\eqref{eq:inertia-clique-cover}. Bollob\'as's chromatic-number
theorem~\cite{Bollobas} gives $\chi(\overline G)=O(n/\log n)$ with
high probability. In Subsection~\ref{sec:matching-upper}, we combine
a standard quantitative form of this coloring bound with
Theorem~\ref{thm:main-random}(\ref{item:random-lower}) to obtain matching
bounds for $n_{\geq0}(G)$ with an explicit failure probability.

More precisely, there are absolute constants $c,c_0>0$ and a
deterministic sequence $\varepsilon_n\geq0$ with $\varepsilon_n\to0$ such that, for all
sufficiently large $n$, with probability at least
$1-\exp(-cn\log(2n))$, a graph $G\sim G(n,1/2)$ satisfies
\begin{equation}\label{eq:random-order}
 c_0\frac{n}{\log(2n)}
 \leq n_{\geq0}(G)
 \leq\left(\frac{\log 2}{2}+\varepsilon_n\right)
       \frac{n}{\log n}.
\end{equation}
In particular, $n_{\geq0}(G)=\Theta(n/\log n)$ on this event.

For $G\sim G(n,1/2)$, we can now compare the optimal inertia bound
with $\alpha(G)$. The clique-number estimate of Bollob\'as and
Erd\H{o}s~\cite{BollobasErdos}, applied to $\overline G$, gives
\begin{equation}\label{eq:random-independence}
 \alpha(G)=2\log_2 n-2\log_2\log_2 n+O(1)
 \qquad\text{with high probability}.
\end{equation}

Combining this known independence-number estimate with
\eqref{eq:random-order} proves Theorem~\ref{thm:main-random}(\ref{item:random-gap}).
It also expresses the optimal inertia bound as a function of
$\alpha(G)$: with high probability,
\begin{equation}\label{eq:random-inertia-alpha}
 n_{\geq0}(G)=\Theta\!\left(2^{\alpha(G)/2}\right).
\end{equation}

\begin{proof}[Proof of Theorem~\ref{thm:main-random}(\ref{item:random-gap})]
The events in~\eqref{eq:random-order}
and~\eqref{eq:random-independence} hold together with probability
$1-o(1)$. On their intersection, dividing $n_{\geq0}(G)$ by
$\alpha(G)=\Theta(\log n)$ gives a ratio of order $n/(\log n)^2$.
To obtain~\eqref{eq:random-inertia-alpha}, observe that
$2^{\alpha(G)/2}=2^{\log_2 n-\log_2\log_2 n+O(1)}
=\Theta(n/\log n)$.
Combining this with~\eqref{eq:random-order} gives the claimed estimate.
\end{proof}

Theorem~\ref{thm:main-random}(\ref{item:random-gap}) gives a quantitative
negative answer to Ihringer's question. Although each prescribed
matrix gives an efficiently evaluable inertia bound, even the
smallest such bound typically exceeds $\alpha(G)$ by a factor of order
$n/(\log n)^2$. Equivalently, the optimal inertia bound grows
exponentially in the independence number on this family of graphs.

The conclusion of~\eqref{eq:random-order} concerning the order of
the inertia bound extends to every fixed $p\in(0,1)$.
For $G\sim G(n,p)$ and all sufficiently large $n$, $n_{\geq0}(G)$ is of order
$n/\log n$ with probability at least $1-\exp(-\delta_p n\log(2n))$.
Here $\delta_p>0$ and the implicit constants may depend on $p$.
The proof is given in Subsection~\ref{sec:matching-upper}.

All inertia lower bounds in this paper also apply when the minimum
in~\eqref{eq:inertia-definition} is restricted to real symmetric matrices.
Since the clique-partition construction uses real symmetric matrices,
the random-graph conclusions and equality for perfect graphs remain
valid under this restriction.

\subsection{Large \texorpdfstring{$n_{\geq0}(G)$}{n>=0(G)} for graphs with small \texorpdfstring{$\alpha(G)$}{alpha(G)}}

The random-graph result concerns graphs whose independence number
grows with their order. Our second objective is to quantify the
separation when the independence number stays fixed. The result of
Tang, Zhang, and Elphick cited above shows that
$n_{\geq0}(G)$ cannot be bounded by any function of $\alpha(G)$ alone;
we strengthen the quantitative separation for $\alpha(G)=2$.

Section~\ref{sec:concluding} develops deterministic lower bounds
for $n_{\geq0}(G)$ from spectral information about the complement.
Applying them to Alon's construction~\cite{Alon1994} yields the
following theorem, improving the exponent $1/9$
in~\cite[Theorem~1.2]{TZE}.

\begin{theorem}\label{thm:alpha-two}
For every sufficiently large $n$, there is an $n$-vertex
graph $G_n$ such that
\[
 \alpha(G_n)=2,\qquad
 n_{\geq0}(G_n)=\Omega(n^{1/3}).
\]
\end{theorem}

The proof is given in Subsection~\ref{sec:applications}.
Applying the same method to Alon's graph family~\cite{Alon2019}
yields the following theorem for graphs with Shannon capacity
$\Theta(G)=3$, and hence $\alpha(G)\leq3$.
The definition of Shannon capacity is recalled in
Subsubsection~\ref{sec:capacity}.

\begin{theorem}\label{thm:capacity}
For every power of two $q\geq4$, there is a graph $H_q$
of order $n=q^2-1$ such that
\[
 \Theta(H_q)=3
\]
and
\begin{equation}\label{eq:capacity-inertia}
 n_{\geq0}(H_q)\geq
 \frac{\sqrt q+1}{2}
 =\Omega(n^{1/4}).
\end{equation}
\end{theorem}

The proof is also given in Subsection~\ref{sec:applications}.
Thus even when the largest independent set has bounded size,
every admissible inertia bound can grow polynomially with the
number of vertices.

The remainder of the paper is organized as follows.
Section~\ref{sec:preliminaries} establishes the normalization
lemma and the coordinate estimates.
Section~\ref{sec:random-energy} proves
Theorem~\ref{thm:main-random}(\ref{item:random-lower}) and establishes
the matching upper bound in~\eqref{eq:random-order}.
The deduction of Theorem~\ref{thm:main-random}(\ref{item:random-gap})
from these bounds is given above.
Section~\ref{sec:concluding} develops the deterministic
complement spectral bounds and their applications to graphs with
bounded independence number.
\section{Preliminaries}\label{sec:preliminaries}

\paragraph{Notation.}
The matrices $I$ and $J$ are the identity and all-ones matrices
of the indicated order, and $\one$ is the all-ones vector.
We write $\rho(M)$ for the spectral radius of a square matrix $M$,
that is, the largest absolute value of its eigenvalues.
We use $\|M\|_F^2=\sum_{i,j}|M_{ij}|^2$ for its squared Frobenius norm,
and $M[S]$ for the principal submatrix indexed by $S$.
The notation $M\succeq0$ means that $M$ is Hermitian positive
semidefinite.

\subsection{Diagonal normalization and Schur-product comparison}\label{sec:normalization}

For a Hermitian matrix $M$, write $|M|=(M^2)^{1/2}$ and let
$M_+$ and $M_-$ denote its positive and negative spectral parts.
Thus $M=M_+-M_-$, $M_+,M_-\succeq0$, $M_+M_-=0$, and
$|M|=M_++M_-$.
The next lemma constructs a matrix $B$ with
$\rank B_+=n_{\geq0}(A)$ and $\tr B_+=n/2$.
The random-support estimate will be applied to the scaled matrix
$DAD$, which has the same zero pattern as $A$.

\begin{lemma}\label{lem:balance}
Let $A$ be an $n\times n$ zero-diagonal Hermitian matrix, where $n\geq1$.
There are $\varepsilon>0$ and a positive diagonal matrix
$D=\diag(d_1,\ldots,d_n)$ such that, writing
$B\triangleq D(A+\varepsilon I)D=B_+-B_-$,
we have $\rank B_+=n_{\geq0}(A)$ and
\begin{equation}\label{eq:balance}
 (B_+)_{ii}=(B_-+\varepsilon D^2)_{ii}=\frac12
 \qquad(1\leq i\leq n).
\end{equation}
In particular, $DAD$ has exactly the same zero entries as $A$.
\end{lemma}

\begin{proof}
Choose $\varepsilon>0$ smaller than the absolute value of every
negative eigenvalue of $A$; if there are none, any positive
$\varepsilon$ will do. Let $C=A+\varepsilon I$. Then
$n_+(C)=n_{\geq0}(A)$ and $C_{ii}=\varepsilon>0$.
Consider the set
\[
 \mathcal D
 =
 \left\{
 S=\diag(s_1,\ldots,s_n):
 s_i>0,\ \prod_{i=1}^n s_i=1
 \right\},
\]
and define $\Phi(S)=\tr\bigl((SCS)_+\bigr)$ for $S\in\mathcal D$.

\smallskip\noindent
\textbf{Claim 1.}
\emph{The function $\Phi$ attains its minimum on $\mathcal D$.}

For every $S\in\mathcal D$,
$\Phi(S)\geq\tr(SCS)=\varepsilon\sum_{i=1}^n s_i^2$.
Thus, on any nonempty sublevel set of $\Phi$, all $s_i$ are
bounded above. Since $\prod_i s_i=1$, they are also bounded
away from zero. Since $\Phi$ is continuous, its nonempty sublevel
sets are compact, and it attains its minimum at some $S_0\in\mathcal D$.

Fix such a minimizer $S_0=\diag(s_1,\ldots,s_n)$ and let
$H=S_0CS_0$.

\smallskip\noindent
\textbf{Claim 2.}
\emph{The diagonal entries of $H_+$ are all equal and positive.}

To compare the value of $\Phi$ before and after rescaling, we use
the inequality $\tr\bigl((U-V)_+\bigr)\leq\tr U$ for $U,V\succeq0$.
Indeed, if $\Pi$ is the positive spectral projection of $U-V$, then
\[
 \tr\bigl((U-V)_+\bigr)
 =\tr(\Pi U)-\tr(\Pi V)
 \leq\tr(\Pi U)\leq\tr U,
\]
because $0\preceq\Pi\preceq I$ and $U,V\succeq0$.

Write $h_i=(H_+)_{ii}$. Since $H=H_+-H_-$ and $H_-\succeq0$,
$h_i\geq H_{ii}=\varepsilon s_i^2>0$.
Let
\[
 g=\left(\prod_{i=1}^n h_i\right)^{1/n},
 \qquad
 E=\diag\left(
 \sqrt{\frac{g}{h_1}},\ldots,\sqrt{\frac{g}{h_n}}
 \right).
\]
Then $\det E=1$, so $ES_0\in\mathcal D$, and
$(EH_+E)_{ii}=g$ for $1\leq i\leq n$.
Since $E$ and $S_0$ commute, the trace inequality above gives
\[
 \begin{aligned}
 \Phi(ES_0)
 &=\tr\bigl((EHE)_+\bigr)\\
 &=\tr\bigl((EH_+E-EH_-E)_+\bigr)\\
 &\leq\tr(EH_+E)\\
 &=ng\\
 &\leq\sum_{i=1}^n h_i
 =\Phi(S_0),
 \end{aligned}
\]
where the last inequality is the arithmetic--geometric mean
inequality. Since $S_0$ is a minimizer, equality holds throughout.
Equality in the arithmetic--geometric mean inequality therefore gives
$h_1=\cdots=h_n=g>0$, proving Claim 2.

Now set
$D=\frac{S_0}{\sqrt{2g}}$ and $B=DCD=\frac{H}{2g}$.
Positive scalar multiplication commutes with taking spectral parts,
so
$B_+=\frac{H_+}{2g}$ and $B_-=\frac{H_-}{2g}$, giving
$(B_+)_{ii}=\frac12$.
Since $D$ is invertible, Sylvester's law of inertia
(\cite[Lemma~2.4]{KW}) gives
$\rank B_+=n_+(B)=n_+(C)=n_{\geq0}(A)$.

The zero diagonal of $A$ yields
$B_{ii}=\varepsilon d_i^2>0$,
and therefore
$(B_-)_{ii} =(B_+)_{ii}-B_{ii} =\frac12-\varepsilon d_i^2$.
Since $B_-\succeq0$, we obtain
$0<\varepsilon d_i^2\leq\frac12$ and $0\leq(B_-)_{ii}<\frac12$.
Finally,
\[
 DAD=B-\varepsilon D^2
   =B_+-(B_-+\varepsilon D^2),
\]
where $B_-+\varepsilon D^2\succeq0$ and
$(B_-+\varepsilon D^2)_{ii} =(B_-)_{ii}+\varepsilon d_i^2 =\frac12$.
As $(DAD)_{ij}=d_ia_{ij}d_j$ with all $d_i>0$, the zero entries
of $DAD$ and $A$ coincide.
\end{proof}

The spectral decomposition $B=B_+-B_-$ gives $B_+B_-=0$.
Subtracting $\varepsilon D^2$ restores the zero diagonal and yields
$DAD=B_+-(B_-+\varepsilon D^2)$, a difference of two positive
semidefinite matrices each with diagonal entries $1/2$.
In general, $B_+$ and $B_-+\varepsilon D^2$ are not the spectral
parts of $DAD$. For each fixed $A$, the shift $\varepsilon$ is
chosen once and kept fixed in the subsequent argument.

For the complement spectral bounds in Section~\ref{sec:concluding},
we use the following lemma to bound the largest eigenvalue of $B_+$.

\begin{lemma}\label{lem:Schur}
Let $B$ be a Hermitian matrix of order $n\geq1$, with positive
and negative spectral parts $B_+$ and $B_-$. Suppose that
$(B_+)_{ii}=\frac12$ and $B_{ii}\geq0$ for $1\leq i\leq n$.
Let $\circ$ denote entrywise multiplication. Suppose a Hermitian
matrix $F$ satisfies $F\circ B=0$ and
\begin{equation}\label{eq:Schur-sandwich}
 aI+cJ\preceq F\preceq bI+rJ
\end{equation}
for real numbers $a,b,c,r$ with
$c>0$, $b\geq0$, and $a<b$.
Then
\begin{equation}\label{eq:Schur-rank}
 \rho(B_+)\leq\frac{b-a}{2c},
 \qquad
 \rank B_+\geq\frac{nc}{b-a}.
\end{equation}
\end{lemma}

\begin{proof}
Set
$\Delta=\Diag(B_{11},\ldots,B_{nn})$.
By the assumptions,
$\Delta\succeq0$, $I\circ B_+=\frac12 I$, and $I\circ B_-=\frac12 I-\Delta$.
The Schur product theorem~\cite{Schur} says that the entrywise
product of two Hermitian positive semidefinite matrices is
positive semidefinite. Applying it to the lower bound
in~\eqref{eq:Schur-sandwich} with $B_+$ gives
$\frac a2 I+cB_+\preceq F\circ B_+$.
Applying it to the upper bound with $B_-$ gives
$F\circ B_-\preceq b\left(\frac12 I-\Delta\right)+rB_-$.
Since $B=B_+-B_-$ and $F\circ B=0$, we have
$F\circ B_+=F\circ B_-$.
Combining the two inequalities and using $b\geq0$ and $\Delta\succeq0$
gives
\[
 cB_+-rB_-
 \preceq \frac{b-a}{2}I-b\Delta
 \preceq \frac{b-a}{2}I.
\]

Let $u$ be a unit eigenvector corresponding to the largest
eigenvalue of $B_+$. This eigenvalue is positive because
$\tr B_+=n/2$.
The identity $B_-B_+=0$ therefore implies $B_-u=0$.
Evaluating the preceding inequality at $u$ gives
$c\rho(B_+)\leq\frac{b-a}{2}$.
Finally,
$\frac n2=\tr B_+ \leq(\rank B_+)\rho(B_+)$
proves the rank bound.
\end{proof}

\subsection{The factorization norm}\label{sec:factorization}

To apply the coordinate-section estimate to matrices, we represent
their entries as inner products of bounded vectors. We use the
standard factorization norm; see \cite[Section~2, equation~(2)]{SS05}
and \cite[Section~2.1]{LMSS}.

\begin{definition}\label{def:gamma-two}
For a real or complex $r\times c$ matrix $M=(M_{ij})$, define
\[
 \gamma_2(M)\triangleq \inf\left\{
  \max_i\|u_i\|_2\max_j\|v_j\|_2:
  M_{ij}=\langle u_i,v_j\rangle\right\}.
\]
The factorizations use a real Hilbert space for a real matrix and a
complex Hilbert space for a complex matrix; complex inner products
are conjugate-linear in the first argument.
\end{definition}

\begin{lemma}\label{lem:gamma-properties}
Let $M$ be a real or complex $r\times c$ matrix.
\begin{enumerate}
\item For every scalar $\lambda$ over the same field,
$\gamma_2(\lambda M)=|\lambda|\gamma_2(M)$, and
$\gamma_2(M)>0$ whenever $M\ne0$.
\item The condition $\gamma_2(M)\leq1$ is equivalent to a
factorization $M_{ij}=\langle u_i,v_j\rangle$ with
$\|u_i\|_2,\|v_j\|_2\leq1$.
Equivalently, there are matrices $P,Q$ over the same field such that
\[
 X=\begin{pmatrix}P&M\\M^*&Q\end{pmatrix}\succeq0,
 \qquad X_{kk}\leq1\quad(1\leq k\leq r+c).
\]
\end{enumerate}
\end{lemma}

\begin{proof}
Rescaling the factors proves homogeneity, while Cauchy--Schwarz gives
$\max_{i,j}|M_{ij}|\leq\gamma_2(M)$ and hence positivity for $M\ne0$.

To justify the second assertion also when $\gamma_2(M)=1$, we show
that the infimum in Definition~\ref{def:gamma-two} is attained.
For $M\ne0$, rescale the two vector families in a minimizing sequence
by reciprocal positive scalars so that their maximal norms are equal.
Their Gram matrices then have fixed order $r+c$ and bounded entries.
A subsequence converges to a positive semidefinite matrix with cross
block $M$, whose Gram representation attains the infimum.
For $M=0$, take all vectors to be zero
(see also \cite[Section~3.3]{LMSS}).
Thus $\gamma_2(M)\leq1$ gives a factorization with all vector norms
at most one; the converse follows directly from the definition.
Finally, the Gram matrix of $u_1,\ldots,u_r,v_1,\ldots,v_c$ has the
stated block form, and a Gram representation of any such $X$ yields
the required vectors.
\end{proof}

Throughout, a matrix, whether or not it is Hermitian, is
\emph{supported on $G$} if its diagonal entries and all entries
corresponding to nonedges are zero.
\subsection{The log-determinant potential for Gram matrices}
\label{sec:logdet}

For the Gram representations in Lemma~\ref{lem:gamma-properties},
we use the negative log-determinant as a potential whose curvature
controls changes in individual entries. This property holds with a
constant independent of the matrix size, as shown in
\cite[Lemma~8 in the full version]{Christiano} and \cite[Proposition~5.1]{MHT}. We give the proof for the coefficient
$1/32$ in~\eqref{eq:potential-curvature}.

\begin{lemma}\label{lem:logdet}
Let $m\geq1$, and define
\[
 K_m=\{X\in\R^{m\times m}:X\succeq0,\ X_{ii}\leq1\},
 \qquad
 f(X)=-\log\det(I+X).
\]
Then $K_m$ is a nonempty compact convex set, and $f$ is continuous
and strictly convex on $K_m$. Moreover,
\begin{equation}\label{eq:logdet-range}
 -m\log2\leq f(X)\leq0
 \qquad(X\in K_m).
\end{equation}
For any linear map $\pi:\R^{m\times m}\to\R^N$ selecting matrix
entries, we have $\pi(K_m)\subseteq[-1,1]^N$ and
\begin{equation}\label{eq:potential-curvature}
 f\left(\frac{U+V}{2}\right)
 \leq
 \frac{f(U)+f(V)}2
 -\frac1{32}\|\pi(U)-\pi(V)\|_\infty^2
 \qquad(U,V\in K_m).
\end{equation}
\end{lemma}

\begin{proof}
The set $K_m$ contains the zero matrix and is closed and convex.
For every $X\in K_m$, positive semidefiniteness gives
$|X_{ij}|^2\leq X_{ii}X_{jj}\leq1$.
Thus $K_m$ is bounded and hence compact. This also shows that
$\pi(K_m)\subseteq[-1,1]^N$.

Since $I+X$ is positive definite for every $X\in K_m$, the function
$f$ is well defined and continuous on $K_m$. All eigenvalues of
$I+X$ are at least one, while Hadamard's inequality gives
$1\leq\det(I+X) \leq\prod_{i=1}^m(1+X_{ii}) \leq2^m$.
Taking negative logarithms proves~\eqref{eq:logdet-range}.

It remains to prove strict convexity and the stated midpoint
inequality. Take $U,V\in K_m$, and let
$H_0=V-U$, $Y=I+U+sH_0$, $g(s)=f(U+sH_0)$, $0\leq s\leq1$.
Along this segment, $Y$ is positive definite and
$Y_{ii}=1+(1-s)U_{ii}+sV_{ii}\leq2$.
Set
$W=Y^{-1/2}H_0Y^{-1/2}$.
The standard Hessian formula for $-\log\det$
(\cite[Appendix~A.4.3, p.~644]{BV}) gives
\[
 g''(s)
 =\tr(Y^{-1}H_0Y^{-1}H_0)
 =\tr(W^2)
 =\|W\|_F^2,
\]
where the second equality follows from cyclicity of the trace
and the last from the symmetry of $W$.
If $U\ne V$, then $W\ne0$, so $g''(s)>0$.
Hence $f$ is strictly convex on $K_m$.

Let $e_1,\ldots,e_m$ be the standard basis vectors of $\R^m$.
Since $H_0=Y^{1/2}WY^{1/2}$ and
$\|Y^{1/2}e_i\|_2^2=Y_{ii}\leq2$, Cauchy--Schwarz gives
\[
 \begin{aligned}
 |(H_0)_{ij}|
 &=\left|(Y^{1/2}e_i)^{\mathsf T}W(Y^{1/2}e_j)\right|\\
 &\leq\sqrt{Y_{ii}Y_{jj}}\,\|W\|_F
 \leq2\|W\|_F.
 \end{aligned}
\]
Since $\pi$ selects matrix entries, it follows that
\[
 g''(s)=\|W\|_F^2
 \geq\frac14\|\pi(V)-\pi(U)\|_\infty^2.
\]
Write
$\kappa=\frac14\|\pi(V)-\pi(U)\|_\infty^2$.
The lower bound $g''(s)\geq\kappa$ implies that
$g(s)-\kappa s^2/2$ is convex on $[0,1]$. Its midpoint inequality gives
$g(1/2)-\frac{\kappa}{8} \leq \frac{g(0)+g(1)-\kappa/2}{2}$,
and therefore
$g(1/2) \leq\frac{g(0)+g(1)}2-\frac{\kappa}{8}$.
Substituting the definitions of $g$ and $\kappa$ gives
\[
 f\left(\frac{U+V}{2}\right)
 \leq
 \frac{f(U)+f(V)}2
 -\frac1{32}\|\pi(U)-\pi(V)\|_\infty^2,
\]
as required.
\end{proof}

\subsection{Coordinate cuts on a random graph}

We now apply the log-determinant potential to Gram matrices whose
cross blocks are supported on a random graph. We identify an
$n\times n$ real matrix with a vector in $\R^{n^2}$ by listing its
entries in a fixed order. The estimate below holds simultaneously
for all such cross blocks, including those chosen after the graph
is revealed.

\begin{lemma}\label{lem:potential-section}
Let $n\geq2$ and $p\in(0,1)$, and set $N=n^2$,
\[
 K=K_{2n}=\{X\in\R^{2n\times2n}:X\succeq0,\ X_{kk}\leq1\},
 \qquad f(X)=-\log\det(I+X),
\]
and
\[
 \pi(X)=(X_{i,n+j})_{i,j\in[n]}\in\R^N.
\]
There is a constant $C_p>0$, depending only on $p$, such that, for
every $t\geq0$, with probability at least $1-e^{-t}$, a graph
$G\sim G(n,p)$ satisfies
\begin{equation}\label{eq:potential-energy}
 \|\pi(X)\|_2^2
 =\sum_{i,j=1}^n|X_{i,n+j}|^2
 \leq C_p\bigl[n\log(2n)+t\bigr]
\end{equation}
for every $X\in K$ whose cross block is supported on $G$.
\end{lemma}

\begin{proof}
Let $\Omega\subseteq[N]$ consist of the diagonal positions and both
orientations of every nonedge of $G$. Thus a cross block is supported
on $G$ precisely when its coordinate vector is zero on $\Omega$.
Set $a=\tfrac12\log(1/p)>0$. For a fixed set $S$ of off-diagonal
positions, let $e(S)$ be the number of unordered pairs it represents.
Since each pair accounts for at most two positions, $e(S)\geq|S|/2$,
and independence of the edge indicators gives
$\mathbb P(S\cap\Omega=\varnothing)=p^{e(S)}\leq e^{-a|S|}$.
If $S$ contains a diagonal position, this probability is zero.
Consequently, for every fixed $S\subseteq[N]$,
\begin{equation}\label{eq:mask-assumption}
 \mathbb P(S\cap\Omega=\varnothing)\leq e^{-a|S|}.
\end{equation}

By Lemma~\ref{lem:logdet}, $K$ is compact and convex, $f$ is continuous
and strictly convex, and $\pi(K)\subseteq[-1,1]^N$.
Moreover, $\max_K f-\min_K f\leq R=2n\log2$, and
\eqref{eq:potential-curvature} holds for this $f$ and $\pi$.
The set $K$, the function $f$, and the map $\pi$ are fixed before
$G$ is sampled.

For each nonempty compact convex set $C\subseteq K$, let $U(C)$
be the unique minimizer of $f$ on $C$, and let
$h(C)=\pi(U(C))$.
Existence follows from compactness and continuity, and uniqueness
from strict convexity. We use $h(C)$ as the representative of $C$.

We shall form nested sets by cutting off the minimizer in one
coordinate at a time. To see how much a cut increases the minimum
of $f$, suppose that
\[
 C'=
 \left\{
 V\in C:
 s\bigl(\pi_i(V)-h_i(C)\bigr)\geq\eta
 \right\}
\]
is nonempty, where $i\in[N]$, $s\in\{-1,1\}$, and $\eta>0$.
Write $U=U(C)$ and $U'=U(C')$. The cut separates $\pi(U')$ from
$\pi(U)$ in coordinate $i$, so $\|\pi(U')-\pi(U)\|_\infty\geq\eta$.
Since $C$ is convex, the midpoint of $U$ and $U'$ belongs to $C$.
The minimality of $U$ on $C$ and~\eqref{eq:potential-curvature} give
\[
 f(U)
 \leq f\!\left(\frac{U+U'}2\right)
 \leq \frac{f(U)+f(U')}{2}-\frac{\eta^2}{32}.
\]
Consequently,
\[
 f(U')-f(U)\geq\frac{\eta^2}{16}.
\]

For a particular target $x$, the useful cuts will depend on $x$.
To apply~\eqref{eq:mask-assumption} before choosing $x$, we first
construct a tree containing all the cuts we will use. Its root is $K$.
At a node $C$, consider all cuts
\[
 C'=
 \left\{
 V\in C:
 s\bigl(\pi_i(V)-h_i(C)\bigr)
 \geq\frac{j}{2\sqrt N}
 \right\},
 \qquad i,j\in[N],\quad s\in\{-1,1\},
\]
and retain each nonempty set as a child. Each child is compact and
convex. There are at most $2N^2$ choices at each node, so the number of representatives
at depth $q$ is at most $(2N^2)^q$. The tree and all its
representatives are fixed before $\Omega$ is sampled.

We next choose how many coordinates of each representative to inspect.
Their number will increase with the depth to offset the growth in
the number of representatives. For a fixed $t\geq0$, set
\[
 b=\left\lceil\frac{t+\log2}{a}\right\rceil,
 \qquad
 d=\left\lceil\frac{2\log(2N)}a\right\rceil.
\]
Then $e^{-ab}\leq\frac{e^{-t}}2$ and $2N^2e^{-ad}\leq\frac12$.
For each representative $h$ at depth $q$ with $b+dq<N$,
let $S(h,q)$ consist of the $b+dq+1$ indices with largest
$|h_i|$, breaking ties by index. Since $S(h,q)$ is fixed,
\eqref{eq:mask-assumption} bounds the probability that it misses
$\Omega$. Taking a union bound over the representatives at all
depths, we find that the probability that any of these sets misses
$\Omega$ is at most
\[
 \begin{aligned}
 \sum_{q\geq0}(2N^2)^q e^{-a(b+dq+1)}
 &=e^{-a(b+1)}
   \sum_{q\geq0}\bigl(2N^2e^{-ad}\bigr)^q\\
 &\leq2e^{-a(b+1)}
 \leq e^{-t}.
 \end{aligned}
\]
Depths with $b+dq\geq N$ contribute no selected sets;
including them in the sum only enlarges the bound.

For the rest of the proof, fix a realization of $\Omega$ for which
every selected set $S(h,q)$ meets $\Omega$. This event has probability
at least $1-e^{-t}$ and was defined without reference to a target
vector.

Consider any vector $x$ that is zero on $\Omega$. If at least
$b+dq+1$ of its coordinates have absolute value at least $r>0$,
then every representative $h$ at depth $q$ satisfies
\[
 \|x-h\|_\infty\geq\frac r2.
\]
Indeed, if $\|x-h\|_\infty<r/2$, then at least $b+dq+1$
coordinates of $h$ have absolute value greater than $r/2$.
Hence $|h_i|>r/2$ for every $i\in S(h,q)$. But for
$i\in S(h,q)\cap\Omega$ we have $x_i=0$, giving
$|h_i|=|h_i-x_i|<\frac r2$, a contradiction.

Now fix any $X^*\in K$ whose cross block is supported on $G$, and
set $x=\pi(X^*)$. Then $x_i=0$ for all $i\in\Omega$, and $X^*$ may
be chosen after $G$ is revealed. Let
$\rho_1\geq\rho_2\geq\cdots\geq\rho_N$
be the absolute values of the coordinates of $x$, arranged
in decreasing order. After the first $b$ terms, we shall use each
successive cut to estimate the sum of the squares in a block of at
most $d$ terms.

The first $\min(b,N)$ terms contribute at most $b$, since
$\rho_k\leq1$. If $b\geq N$, this bounds the entire sum and gives
the claim by the choice of $b$. Otherwise, we follow a path from
$C_0=K$ in the fixed tree, choosing at each step a child that still
contains $X^*$. At depth $q$, write
$U_q=U(C_q)$, $h^{(q)}=\pi(U_q)$.

At depth $q$, stop if $b+dq\geq N$, since all terms have been
counted. Otherwise the largest remaining term is $r_q=\rho_{b+dq+1}$.
If $r_q\leq N^{-1/2}$, the remaining sum of the squares is at most
$Nr_q^2\leq1$, and we stop.

Suppose, therefore, that $r_q>N^{-1/2}$. There are at least
$b+dq+1$ coordinates of $x$ with absolute value at least $r_q$.
The estimate for representatives at depth $q$, applied to $h^{(q)}$,
gives a coordinate $i\in[N]$ and a sign $s\in\{-1,1\}$ such that
$s\bigl(x_i-h_i^{(q)}\bigr)\geq\frac{r_q}{2}$.
To choose a cut already present in the tree, take
$j=\lfloor r_q\sqrt N\rfloor$ and $\eta=\frac{j}{2\sqrt N}$.
Since $N^{-1/2}<r_q\leq1$, we have
$1\leq j\leq N$ and $\frac{r_q}{4}\leq\eta\leq\frac{r_q}{2}$.
Take
\[
 C_{q+1}
 =
 \left\{
 V\in C_q:
 s\bigl(\pi_i(V)-h_i^{(q)}\bigr)\geq\eta
 \right\}.
\]
Since $s\bigl(\pi_i(X^*)-h_i^{(q)}\bigr) =s\bigl(x_i-h_i^{(q)}\bigr) \geq\frac{r_q}{2}\geq\eta$,
this cut contains $X^*$. Thus $C_{q+1}$ is a nonempty child of $C_q$.

Applying the increment estimate for a cut of size $\eta$ gives
\[
 f(U_{q+1})-f(U_q)
 \geq\frac{\eta^2}{16}
 \geq\frac{r_q^2}{256}.
\]
The next block has at most $d$ terms, none larger than $r_q$.
The increase in the minimum of $f$ therefore bounds the contribution
of this block to the sum of the squares, giving
\[
 \begin{aligned}
 \sum_{k=b+dq+1}^{\min\{N,b+d(q+1)\}}\rho_k^2
 &\leq d r_q^2\\
 &\leq256d\bigl(f(U_{q+1})-f(U_q)\bigr).
 \end{aligned}
\]

The construction stops after at most
$\lceil (N-b)/d\rceil$ cuts. Let $T$ denote the number of cuts.
Since $X^*\in C_T$ and $U_T$ minimizes $f$ on $C_T$,
$f(U_T)\leq f(X^*)$.
Also $f(U_0)=\min_K f$, and hence
$f(U_T)-f(U_0) \leq f(X^*)-\min_K f \leq R$.
The increases in the minimum of $f$ telescope. Adding the first
$\min(b,N)$ terms and the possible final remainder gives
\[
 \begin{aligned}
 \|x\|_2^2
 &\leq b+
 256d\sum_{q=0}^{T-1}
 \bigl(f(U_{q+1})-f(U_q)\bigr)+1\\
 &=b+256d\bigl(f(U_T)-f(U_0)\bigr)+1\\
 &\leq b+256Rd+1.
 \end{aligned}
\]
Finally,
$b\leq1+\frac{t+\log2}{a}$, $d\leq1+\frac{2\log(2N)}a$.
Substituting $R=2n\log2$, $N=n^2$, and $a=\tfrac12\log(1/p)$ gives
$\|\pi(X^*)\|_2^2\leq C_p\bigl[n\log(2n)+t\bigr]$
for a sufficiently large constant $C_p$ depending only on $p$.
The path through the tree was chosen for $X^*$, but the event is the
same for every $X\in K$ whose cross block is supported on $G$.
\end{proof}

\section{Proof of Theorem~\texorpdfstring{\ref{thm:main-random}(\ref{item:random-lower})}{\getrefnumber{thm:main-random}(\getrefnumber{item:random-lower})}}
\label{sec:random-energy}

We use the log-determinant potential to bound the Frobenius norm of
every matrix supported on $G$ in terms of its factorization norm.
Diagonal normalization then converts this estimate into an inertia bound.

\subsection{Uniform energy on a random graph support}

\begin{theorem}
\label{thm:gamma-energy}
There are absolute constants $C,c>0$ with the following property.
For every $n\geq2$, let $\mathcal E_n$ be the set of graphs $G$
on $[n]$ for which
\begin{equation}\label{eq:gamma-energy}
 \sum_{i,j}|B_{ij}|^2
 \leq C\gamma_2(B)^2 n\log(2n)
\end{equation}
holds for every complex matrix $B$ supported on $G$.
Then, for $G\sim G(n,1/2)$,
\[
 \mathbb P(G\in\mathcal E_n)
 \geq1-\exp\bigl(-cn\log(2n)\bigr).
\]
\end{theorem}
\begin{proof}
Apply Lemma~\ref{lem:potential-section} with $p=1/2$ and
$t=n\log(2n)$. By Lemma~\ref{lem:gamma-properties}, every real matrix
$B$ with $\gamma_2(B)\leq1$ is the cross block of a matrix in $K_{2n}$.
Hence, on an event of probability at least $1-\exp(-n\log(2n))$,
\[
 \|B\|_F^2\leq Cn\log(2n)
\]
for every such real $B$ supported on $G$, where $C$ is absolute.
The event in Lemma~\ref{lem:potential-section} is common to all cross
blocks, so $B$ and its Gram representation may both be chosen after
$G$ is revealed.

For a complex matrix $B$ supported on $G$ with $\gamma_2(B)\leq1$,
use Lemma~\ref{lem:gamma-properties} to choose factor vectors of norms
at most one and write
$u_i=a_i+ib_i$ and $v_j=c_j+id_j$. The real vectors
$(a_i,b_i)$ and $(c_j,d_j)$ have the same norms as $u_i$ and $v_j$,
and their inner product is $(\operatorname{Re}B)_{ij}$. Thus
$\gamma_2(\operatorname{Re}B)\leq1$. Applying the same argument to $-iB$
also gives $\gamma_2(\operatorname{Im}B)\leq1$.
Both real matrices are supported on $G$, so the preceding estimate
applies to them on the same event. Since
$\|B\|_F^2=\|\operatorname{Re}B\|_F^2+\|\operatorname{Im}B\|_F^2$,
adding these bounds yields $\|B\|_F^2\leq Cn\log(2n)$ after adjusting $C$.
For $B\ne0$, Lemma~\ref{lem:gamma-properties} allows us to apply this
estimate to $B/\gamma_2(B)$ and rescale to obtain~\eqref{eq:gamma-energy}.
The inequality also holds for $B=0$.
\end{proof}

\subsection{From matrix energy to inertia}

Lemma~\ref{lem:balance} associates two matrices to a weight matrix $A$.
The positive part of $B=D(A+\varepsilon I)D$ has rank $n_{\geq0}(A)$,
while $DAD$ has the zero entries needed to apply the energy bound.
\begin{proof}[Proof of Theorem~\ref{thm:main-random}(\ref{item:random-lower})]
Work on the event of Theorem~\ref{thm:gamma-energy}, and fix any
$A\in\W(G)$. Use Lemma~\ref{lem:balance} to write
\[
 \begin{gathered}
 B=B_+-B_-=DAD+\varepsilon D^2,\\
 DAD=B_+-(B_-+\varepsilon D^2),\\
 (B_+)_{ii}=(B_-+\varepsilon D^2)_{ii}=1/2,\qquad
 \rank B_+=n_{\geq0}(A).
 \end{gathered}
\]
Here $B=D(A+\varepsilon I)D$, with $\varepsilon>0$ and $D$ positive
diagonal. The spectral parts $B_+,B_-$ are positive semidefinite,
as is $B_-+\varepsilon D^2$.
Take Gram representations $(B_+)_{ij}=\langle u_i,u_j\rangle$ and
$(B_-+\varepsilon D^2)_{ij}=\langle v_i,v_j\rangle$.
The diagonal identities give $\|u_i\|_2^2=\|v_i\|_2^2=1/2$, and hence
\[
 (DAD)_{ij}=\langle(u_i,v_i),(u_j,-v_j)\rangle,
 \qquad \|(u_i,v_i)\|_2=\|(u_i,-v_i)\|_2=1.
\]
Lemma~\ref{lem:gamma-properties} therefore gives $\gamma_2(DAD)\leq1$.
Diagonal multiplication preserves the zero entries of $A$, so $DAD$ is
supported on $G$ and
Theorem~\ref{thm:gamma-energy} gives $\|DAD\|_F^2\leq S=Cn\log(2n)$.

The spectral parts satisfy $B_+B_-=B_-B_+=0$, so $B^2=B_+^2+B_-^2$.
The zero diagonal of $DAD$ gives
$\tr\bigl((DAD)D^2\bigr)=\tr\bigl(D^2(DAD)\bigr)=0$.
Using that $DAD$ is Hermitian and $0<\varepsilon d_i^2\leq1/2$, we obtain
\[
 \begin{aligned}
 \tr(B_+^2)
 &\leq\tr B^2\\
 &=\tr\bigl((DAD)^2\bigr)+\varepsilon^2\tr(D^4)\\
 &=\|DAD\|_F^2+\varepsilon^2\tr(D^4)
 \leq S+n/4.
 \end{aligned}
\]
Summing the diagonal entries of $B_+$ gives $\tr B_+=n/2$.
If $k=\rank B_+$ and $\mu_1,\ldots,\mu_k$ are its positive
eigenvalues, Cauchy--Schwarz gives
$(\sum_{\ell=1}^k\mu_\ell)^2\leq k\sum_{\ell=1}^k\mu_\ell^2$.
Since $\tr B_+>0$, we obtain
\[
 n_{\geq0}(A)=\rank B_+
 \geq\frac{(\tr B_+)^2}{\tr(B_+^2)}
 \geq\frac{n^2}{4S+n}.
\]
Substituting $S=Cn\log(2n)$ gives the lower bound $n/(4C\log(2n)+1)$.
For $n\geq2$ we have $\log(2n)\geq1$, so this is at least
$(4C+1)^{-1}n/\log(2n)$.
Since $A\in\W(G)$ was arbitrary, the bound holds for all such $A$ on
the event of Theorem~\ref{thm:gamma-energy}.
\end{proof}

\paragraph{Explicit constants.}
The preceding proof permits $c_1=1$ and $c_2=1/32801$ in
Theorem~\ref{thm:main-random}(\ref{item:random-lower}). Indeed, retain
the estimate $b+256Rd+1$ from the proof of Lemma~\ref{lem:potential-section}
and substitute $a=(\log2)/2$, $N=n^2$, $R=2n\log2$, and
$t=n\log(2n)$. On the same event of probability at least
$1-\exp(-n\log(2n))$, every real matrix $B$ supported on $G$
with $\gamma_2(B)\leq1$ satisfies
\[
 \begin{aligned}
 \|B\|_F^2
 &\leq4+\left(4096+\frac2{\log2}\right)n\log(2n)
       -1536n\log2\\
 &\leq\left(4096+\frac2{\log2}\right)n\log(2n).
 \end{aligned}
\]
Adding the estimates for the real and imaginary parts gives an
energy coefficient $8192+4/\log2<8200$. Thus one may take
$C=8200$ and $c=1$ in Theorem~\ref{thm:gamma-energy}, and the
rank estimate above gives $c_2=(4C+1)^{-1}=1/32801>10^{-5}$.

\subsection{Consequences of the lower bound}\label{sec:matching-upper}
We combine Theorem~\ref{thm:main-random}(\ref{item:random-lower}) with
the known upper bound from clique partitions. The probability estimate we need is a
quantitative consequence of Krivelevich--Sudakov--Vu--Wormald's
independent-set bound and the standard coloring argument
in~\cite[Theorem~2.1 and Section~4]{KSVW}: for each fixed $q\in(0,1)$,
there are $b_q>0$ and a deterministic sequence $\eta_n(q)\geq0$
tending to zero such that, for $H\sim G(n,q)$ and all sufficiently
large $n$,
\begin{equation}\label{eq:coloring-tail}
 \mathbb P\!\left(
 \chi(H)>\left(\frac{\log(1/(1-q))}{2}+\eta_n(q)\right)
              \frac{n}{\log n}\right)
 \leq\exp\!\left(-b_q\frac{n^2}{(\log n)^8}\right).
\end{equation}
To verify this estimate, put $m=\lceil n/(\log n)^2\rceil$ and let
\[
 k=k_0(m,q)
 =\max\left\{j:\binom mj(1-q)^{\binom j2}\geq m^4\right\}.
\]
By~\cite[Theorem~2.1]{KSVW},
$k\sim2\log m/\log(1/(1-q))$ and
$\mathbb P(\alpha(G(m,q))<k)\leq\exp(-a_qm^2/(qk^4))$
for some $a_q>0$. A union bound over all $m$-vertex subsets gives
\[
 \begin{aligned}
 &\mathbb P\bigl(\exists S\subseteq V(H),\
     |S|=m,\ \alpha(H[S])<k\bigr)\\
 &\qquad\leq\binom nm
       \exp\!\left(-a_q\frac{m^2}{qk^4}\right)
 \leq\exp\!\left(-b_q\frac{n^2}{(\log n)^8}\right)
 \end{aligned}
\]
for some $b_q>0$ and all sufficiently large $n$, since
$\binom nm\leq2^n$, $k=\Theta_q(\log n)$, and
$n^2/(\log n)^8\gg n$.
On the complementary event, repeatedly remove an independent set
of size $k$ while at least $m$ vertices remain, assigning a new
color to each removed set and then to each remaining vertex.
This gives
\[
 \chi(H)\leq\frac nk+m
 =\left(\frac{\log(1/(1-q))}{2}+o(1)\right)\frac n{\log n}.
\]
The deterministic error can be bounded above by a nonnegative
sequence $\eta_n(q)\to0$, proving~\eqref{eq:coloring-tail}.

\begin{proof}[Proof of~\eqref{eq:random-order}]
Apply~\eqref{eq:coloring-tail} with $q=1/2$ to $\overline G$.
Together with~\eqref{eq:inertia-clique-cover}, this gives the upper
bound in~\eqref{eq:random-order} with $\varepsilon_n=\eta_n(1/2)$.
The lower bound is Theorem~\ref{thm:main-random}(\ref{item:random-lower}),
with $c_0=c_2$.
The two bounds hold simultaneously with probability at least
\[
 1-\exp(-c_1n\log(2n))
   -\exp\!\left(-b_{1/2}\frac{n^2}{(\log n)^8}\right)
 \geq1-\exp(-cn\log(2n))
\]
for some absolute $c>0$ and all sufficiently large $n$, since
$n^2/(\log n)^8\gg n\log(2n)$.
\end{proof}

\begin{remark}
For every fixed $p\in(0,1)$ and $G\sim G(n,p)$,
$n_{\geq0}(G)=\Theta_p(n/\log n)$ with probability at least
$1-\exp(-\delta_p n\log(2n))$ for some $\delta_p>0$ and all
sufficiently large $n$.
This follows by applying Lemma~\ref{lem:potential-section} with edge
probability $p$ and repeating the energy and rank arguments for the
lower bound, and using
\eqref{eq:coloring-tail} with $q=1-p$ together
with~\eqref{eq:inertia-clique-cover} for the upper bound.
\end{remark}

\section{Complement spectral bounds and applications}
\label{sec:concluding}

We use Lemmas~\ref{lem:balance} and~\ref{lem:Schur} to obtain
deterministic complement spectral bounds for $n_{\geq0}(G)$.
We then apply these bounds to two graph families of Alon.

\subsection{Complement spectral bounds}\label{sec:spectral}

\begin{theorem}
\label{thm:regular}
Let $G$ be a graph on $n$ vertices. Suppose that $\overline G$
contains a spanning $d$-regular graph $\Gamma$, where
$0<d<n-1$. Write
$d=\lambda_1\geq\lambda_2\geq\cdots\geq\lambda_n$
for the eigenvalues of $A_\Gamma$. Then
$\lambda_2>\lambda_n$ and
\begin{equation}\label{eq:regular-sharp}
 n_{\geq0}(G)\geq
 \frac{d-\lambda_n}{\lambda_2-\lambda_n}.
\end{equation}
In particular, if $\lambda>0$ satisfies
$|\lambda_i|\leq\lambda$ for every $i\geq2$, then
\begin{equation}\label{eq:regular-simple}
 n_{\geq0}(G)\geq\frac{d+\lambda}{2\lambda}.
\end{equation}
\end{theorem}

Tang--Zhang--Elphick~\cite[Theorem~1.4]{TZE} proved
$n_{\geq0}(G)\geq \left(\frac{d+\lambda}{8\lambda}\right)^{1/3}$
under the additional assumption $d<n/2-1$.
Theorem~\ref{thm:regular} replaces this cube-root dependence by
linear dependence on $d/\lambda$ and removes that restriction
on the degree.
The sharper estimate~\eqref{eq:regular-sharp} uses $\lambda_2$
and $\lambda_n$ separately. Neither estimate requires $\Gamma$
to be connected.

Lemma~\ref{lem:balance} fixes $\tr B_+=n/2$, so an upper bound
on the largest eigenvalue of $B_+$ gives a lower bound on its rank:
$\rank B_+\geq\frac{\tr B_+}{\rho(B_+)}$.
We obtain the required upper bound from Lemma~\ref{lem:Schur},
using the lower spectral bound on $A_\Gamma$ with $B_+$ and
the upper spectral bound with $B_-$.

\begin{proof}[Proof of Theorem~\ref{thm:regular}]
Fix any $A\in\W(G)$.
By Lemma~\ref{lem:balance}, there is a normalized matrix
$B=D(A+\varepsilon I)D$
satisfying~\eqref{eq:balance} and
$\rank B_+=n_{\geq0}(A)$.
Let $F=A_\Gamma$. Its off-diagonal support is disjoint from that
of $B$ because $\Gamma$ is a subgraph of $\overline G$.
Since $F$ also has zero diagonal, $F\circ B=0$.

Regularity gives $F\one=d\one$.
On $\one^\perp$, the eigenvalues of $F$ lie between
$\lambda_n$ and $\lambda_2$. Consequently,
\begin{equation}\label{eq:regular-sandwich}
 \lambda_n I+\frac{d-\lambda_n}{n}J
 \preceq F\preceq
 \lambda_2 I+\frac{d-\lambda_2}{n}J.
\end{equation}
Both bounding matrices agree with $F$ on $\one$, since all three
send $\one$ to $d\one$. On $\one^\perp$ the two multiples of $J$
act as zero, leaving the bounds $\lambda_n I$ and $\lambda_2 I$.

We check the sign conditions needed in Lemma~\ref{lem:Schur}.
Since $d<n-1$, there are two nonadjacent vertices $i,j$ in
$\Gamma$. The vector $e_i-e_j$ belongs to $\one^\perp$
and has Rayleigh quotient zero, giving $\lambda_2\geq0$.
Also, $d>0$ and $\tr F=0$ imply $\lambda_n<0$.
In particular, $\lambda_2-\lambda_n>0$.

We can therefore apply Lemma~\ref{lem:Schur}
to~\eqref{eq:regular-sandwich} with $a=\lambda_n$, $b=\lambda_2$,
$c=\frac{d-\lambda_n}{n}$, and $r=\frac{d-\lambda_2}{n}$. This gives
$n_{\geq0}(A)=\rank B_+ \geq\frac{d-\lambda_n}{\lambda_2-\lambda_n}$.
Taking the minimum over $A\in\W(G)$ proves
\eqref{eq:regular-sharp}.

For~\eqref{eq:regular-simple}, use the bounds $-\lambda$ and
$\lambda$ on the eigenvalues of $F$ on $\one^\perp$. This gives
\[
 -\lambda I+\frac{d+\lambda}{n}J
 \preceq F\preceq
 \lambda I+\frac{d-\lambda}{n}J.
\]
A second application of Lemma~\ref{lem:Schur}, now with
$a=-\lambda$, $b=\lambda$, $c=\frac{d+\lambda}{n}$, $r=\frac{d-\lambda}{n}$,
yields
$n_{\geq0}(A)\geq\frac{d+\lambda}{2\lambda}$.
The lemma imposes no sign condition on $r$, so this argument
also applies when $\lambda>d$.
Taking the minimum over $A$ completes the proof.
\end{proof}

We can also apply the Schur-product comparison to $A_{\overline G}$
itself. For this fixed matrix we choose a centering parameter $p>0$
and bound the spectral radius of $A_{\overline G}-p(J-I)$.

\begin{corollary}
\label{thm:universal}
Let $G$ be a graph on $n\geq2$ vertices.
If $G$ has at least one edge, then for every $p>0$,
\begin{equation}\label{eq:universal-p}
 n_{\geq0}(G)\geq
 \frac{pn}{2\rho\bigl(A_{\overline G}-p(J-I)\bigr)}.
\end{equation}
In particular, the following bound holds for every graph $G$
on $n\geq2$ vertices, including the edgeless graph:
\begin{equation}\label{eq:universal}
 n_{\geq0}(G)\geq
 \frac{n}{2\rho\bigl(2A_{\overline G}-(J-I)\bigr)}.
\end{equation}
\end{corollary}

\begin{proof}
Suppose first that $G$ has at least one edge, and fix $p>0$.
Let $F=A_{\overline G}$, $E=F-p(J-I)$, and $\eta=\rho(E)$.
At an edge $ij$ of $G$ we have $E_{ij}=-p$, so $\eta\geq p>0$.
The inequalities $-\eta I\preceq E\preceq\eta I$ give
\[
 (-p-\eta)I+pJ
 \preceq F\preceq
 (\eta-p)I+pJ.
\]
For any $A\in\W(G)$, Lemma~\ref{lem:balance} gives a matrix
$B$ with $\rank B_+=n_{\geq0}(A)$. Since $F$ is the adjacency
matrix of $\overline G$ and the off-diagonal support of $B$
lies in $G$, we have $F\circ B=0$. Applying Lemma~\ref{lem:Schur}
with $a=-p-\eta$, $b=\eta-p$, and $c=r=p$ yields
$n_{\geq0}(A)=\rank B_+ \geq\frac{pn}{2\eta}$.
Taking the minimum over $A\in\W(G)$ proves
\eqref{eq:universal-p}.

Taking $p=1/2$ in~\eqref{eq:universal-p} yields
\eqref{eq:universal} whenever $G$ has at least one edge.
If $G$ is edgeless, then $n_{\geq0}(G)=n$ and
$A_{\overline G}=J-I$. Since $n\geq2$, we have
\[
 \frac{n}{2\rho\bigl(2A_{\overline G}-(J-I)\bigr)}
 =\frac{n}{2\rho(J-I)}
 =\frac{n}{2(n-1)}
 \leq n
 =n_{\geq0}(G).
\]
Thus~\eqref{eq:universal} holds for every graph on at least
two vertices.
\end{proof}

\subsection{Applications to explicit graphs}\label{sec:applications}

We apply Theorem~\ref{thm:regular} to two graph families of Alon.
In each case we use the known degree and eigenvalue estimates for the
complement to bound the inertia of the original graph.

\subsubsection{Independence number two}

Tang--Zhang--Elphick~\cite[Theorem~1.2]{TZE} proved that
there are $n$-vertex graphs with independence number at most two
and Hermitian inertia $\Omega(n^{1/9})$. Applying Theorem~\ref{thm:regular} to the construction below raises
the exponent to $1/3$.

We recall the construction from~\cite[Section~2]{Alon1994}.
Let $q=2^k$, where $k\geq2$ and $3\nmid k$, and let $n_0=q^3$.
Fix a basis of $\mathbb F_q$ over $\mathbb F_2$, and let
$\ell:\mathbb F_q\to\mathbb F_2$ be the first coordinate map
with respect to this basis. Define
\[
 W_i=\{a\in\mathbb F_q\setminus\{0\}:\ell(a^7)=i\},
 \qquad i\in\{0,1\}.
\]
The graph $\Gamma_{n_0}$ has vertex set
$V(\Gamma_{n_0})=\mathbb F_q^3$.
Two distinct vertices $u,v$ are adjacent if and only if
\[
 u+v=(a+b,\ a^3+b^3,\ a^5+b^5)
\]
for some $a\in W_0$ and $b\in W_1$, where all operations are
performed in $\mathbb F_q$.
At order $n_0$ we take $G_{n_0}=\overline{\Gamma_{n_0}}$.
In the proof we extend this graph to other orders by adding
universal vertices.

\begin{proof}[Proof of Theorem~\ref{thm:alpha-two}]
By Alon~\cite[Theorem~2.1]{Alon1994}, the graph
$\Gamma_{n_0}$ defined above is triangle-free and regular, with degree
\begin{equation}\label{eq:Alon-degree}
 d=2^{k-1}(2^{k-1}-1)=\Theta(n_0^{2/3}),
\end{equation}
and its nonprincipal adjacency eigenvalues have absolute value at
most $\lambda=O(n_0^{1/3})$.
Here $0<d<n_0-1$.
Since $\Gamma_{n_0}$ is triangle-free and has an edge,
$\alpha(\overline{\Gamma_{n_0}})=2$.
Theorem~\ref{thm:regular} gives
\begin{equation}\label{eq:Alon-application}
 n_{\geq0}(\overline{\Gamma_{n_0}})
 \geq\frac{d+\lambda}{2\lambda}
 =\Omega(n_0^{1/3}).
\end{equation}

To pass to arbitrary sufficiently large $n$, we follow
the argument of~\cite[proof of Theorem~1.2]{TZE}.
For $n\geq64$, choose the largest admissible
$n_0=2^{3k}\leq n$.
Consecutive admissible values of $k$ differ by at most two,
so $n_0>n/64$.
Add $n-n_0$ vertices to $\overline{\Gamma_{n_0}}$, joining each
new vertex to every other vertex. Any independent set containing
a new vertex is a singleton, so the resulting graph $G_n$ still
has independence number two. The original graph is an induced
subgraph of $G_n$. Restricting any admissible weight matrix of $G_n$
to this subgraph preserves admissibility and, by eigenvalue
interlacing, cannot increase its nonnegative inertia index. Taking
the minimum over all admissible weight matrices of $G_n$ gives
\[
 n_{\geq0}(G_n)
 \geq n_{\geq0}(\overline{\Gamma_{n_0}})
 =\Omega(n_0^{1/3})
 =\Omega(n^{1/3}).
\]
\end{proof}

Theorem~\ref{thm:alpha-two} gives a lower bound for the
extremal function
$F_2(n)\triangleq \max\{n_{\geq0}(G): |V(G)|=n,\ \alpha(G)=2\}$,
defined for $n\geq2$.
Specifically, $F_2(n)\geq c_3n^{1/3}$ for some absolute
$c_3>0$ and all sufficiently large $n$.
For an upper bound, observe that $\overline G$ is triangle-free
whenever $\alpha(G)=2$. Applying the chromatic-number bound of
Davies and Illingworth~\cite[Theorem~1.1]{DaviesIllingworth}
to $\overline G$ and using~\eqref{eq:inertia-clique-cover},
we obtain
\begin{equation}\label{eq:fixed-independence-range}
 c_3n^{1/3}\leq F_2(n)
 \leq(2+o(1))\sqrt{\frac{n}{\log n}}.
\end{equation}
What is the asymptotic order of $F_2(n)$ within this range?

\subsubsection{Shannon capacity three}\label{sec:capacity}

Recall that the Shannon capacity of a graph $H$ is
$\Theta(H)\triangleq \sup_{t\geq1}\alpha(H^{\boxtimes t})^{1/t}$,
where $H^{\boxtimes t}$ denotes its $t$-fold strong power.
Two distinct tuples are adjacent in this power if their
entries in each coordinate are equal or adjacent in $H$.

Tang--Zhang--Elphick~\cite[Theorem~1.8]{TZE} proved an inertia
lower bound of $\Omega(n^{1/12})$ for a family of graphs with
Shannon capacity at most three. Our bound gives an explicit lower
bound of order $n^{1/4}$ for the same family.

We give an equivalent normalized form of the construction in
\cite[proof of Theorem~3.2]{Alon2019}.
Let $q=2^k\geq4$, and define
\[
 V(H_q)=
 \{(x_0,x_1,x_2)\in\mathbb F_q^3:
   x_0+x_1+x_2=1\}\setminus\{(1,1,1)\}.
\]
For distinct $x,y\in V(H_q)$, set
\[
 x\sim_{H_q}y
 \quad\Longleftrightarrow\quad
 x_0y_0+x_1y_1+x_2y_2\ne0,
\]
where all operations are performed in $\mathbb F_q$.
Thus $|V(H_q)|=q^2-1$, and its complement $\Gamma_q$ has
the same vertex set, with distinct vertices adjacent precisely when
$x_0y_0+x_1y_1+x_2y_2=0$.
To identify this graph with Alon's, observe that each projective class with nonzero coordinate sum has a unique representative
whose coordinate sum is one. In characteristic two the excluded
class of nonzero vectors $(a,a,a)$ is therefore represented by
$(1,1,1)$. Multiplying either representative by a nonzero scalar
preserves orthogonality and hence adjacency.

\begin{proof}[Proof of Theorem~\ref{thm:capacity}]
By Alon~\cite[Theorem~3.2]{Alon2019}, $H_q$ has Shannon capacity
exactly three. Its complement $\Gamma_q$ is $q$-regular, with
$|\lambda_i(A_{\Gamma_q})|\leq\sqrt q$ for $i\geq2$.
These spectral estimates follow from the calculation for the
polarity graph in~\cite[proof of Claim~3.4]{Alon2019} and are
also given explicitly in~\cite[proof of Theorem~1.8]{TZE}.

For $q\geq4$, we have $0<q<q^2-2$.
Thus Theorem~\ref{thm:regular} applies with
$d=q$ and $\lambda=\sqrt q$, giving
$n_{\geq0}(H_q) \geq\frac{q+\sqrt q}{2\sqrt q} =\frac{\sqrt q+1}{2}$.
Since $|V(H_q)|=q^2-1$, this is
$\Omega(|V(H_q)|^{1/4})$.
\end{proof}

\section*{Acknowledgements}

This work was supported in part by the National Natural Science
Foundation of China under Grants 12371354 and W2521102, the
Montenegrin--Chinese Science and Technology Cooperation Project
under Grant 4-3, and the Science and Technology Commission of
Shanghai Municipality under Grant 25LN3200600.


\begin{thebibliography}{99}

\bibitem{Alon1994}
N.~Alon,
\emph{Explicit Ramsey graphs and orthonormal labelings},
Electron. J. Combin. \textbf{1} (1994), Research Paper~12.
\href{https://doi.org/10.37236/1192}{doi:10.37236/1192}.

\bibitem{Alon2019}
N.~Alon,
\emph{Lov\'asz, vectors, graphs and codes},
in \emph{Building Bridges II}, I.~B\'ar\'any, G.~O.~H.~Katona,
and A.~Sali (eds.), Bolyai Society Mathematical Studies~28,
Springer, Berlin--Heidelberg, 2019, pp.~1--16.
\href{https://doi.org/10.1007/978-3-662-59204-5_1}
{doi:10.1007/978-3-662-59204-5\_1}.
Author manuscript:
\url{https://web.math.princeton.edu/~nalon/PDFS/ll70.pdf}.

\bibitem{Bollobas}
B.~Bollob\'as,
\emph{The chromatic number of random graphs},
Combinatorica \textbf{8} (1988), 49--55.
\href{https://doi.org/10.1007/BF02122551}{doi:10.1007/BF02122551}.

\bibitem{BollobasErdos}
B.~Bollob\'as and P.~Erd\H{o}s,
\emph{Cliques in random graphs},
Math. Proc. Cambridge Philos. Soc. \textbf{80} (1976), 419--427.
\href{https://doi.org/10.1017/S0305004100053056}
{doi:10.1017/S0305004100053056}.

\bibitem{BV}
S.~Boyd and L.~Vandenberghe,
\emph{Convex Optimization},
Cambridge University Press, Cambridge, 2004.
\href{https://doi.org/10.1017/CBO9780511804441}
{doi:10.1017/CBO9780511804441}.

\bibitem{Christiano}
P.~Christiano,
\emph{Online local learning via semidefinite programming},
Proceedings of STOC 2014, ACM, 2014, 468--474.
\href{https://doi.org/10.1145/2591796.2591880}
{doi:10.1145/2591796.2591880}.
Full version: arXiv:1403.5287v1.
\url{https://arxiv.org/abs/1403.5287v1}.

\bibitem{DaviesIllingworth}
E.~Davies and F.~Illingworth,
\emph{The $\chi$-Ramsey problem for triangle-free graphs},
SIAM J. Discrete Math. \textbf{36} (2022), 1124--1134.
\href{https://doi.org/10.1137/21M1437573}
{doi:10.1137/21M1437573}.

\bibitem{KSVW}
M.~Krivelevich, B.~Sudakov, V.~H.~Vu, and N.~C.~Wormald,
\emph{On the probability of independent sets in random graphs},
Random Structures Algorithms \textbf{22} (2003), 1--14.
\href{https://doi.org/10.1002/rsa.10063}{doi:10.1002/rsa.10063}.

\bibitem{KW}
M.~Kwan and Y.~Wigderson,
\emph{The inertia bound is far from tight},
Bull. Lond. Math. Soc. \textbf{56} (2024), 3196--3208.
\href{https://doi.org/10.1112/blms.13127}{doi:10.1112/blms.13127}.
\url{https://arxiv.org/abs/2312.04925}.

\bibitem{LMSS}
N.~Linial, S.~Mendelson, G.~Schechtman, and A.~Shraibman,
\emph{Complexity measures of sign matrices},
Combinatorica \textbf{27} (2007), 439--463.
\href{https://doi.org/10.1007/s00493-007-2160-5}
{doi:10.1007/s00493-007-2160-5}.

\bibitem{MHT}
K.-i.~Moridomi, K.~Hatano, and E.~Takimoto,
\emph{Online linear optimization with the log-determinant regularizer},
IEICE Trans. Inf. Syst. \textbf{E101-D} (2018), 1511--1520.
\href{https://doi.org/10.1587/transinf.2017EDP7317}
{doi:10.1587/transinf.2017EDP7317}.
\url{https://arxiv.org/abs/1710.10002}.

\bibitem{Schur}
J.~Schur,
\emph{Bemerkungen zur Theorie der beschr\"ankten Bilinearformen
mit unendlich vielen Ver\"anderlichen},
J. Reine Angew. Math. \textbf{140} (1911), 1--28.
\href{https://doi.org/10.1515/crll.1911.140.1}
{doi:10.1515/crll.1911.140.1}.

\bibitem{Sinkovic}
J.~Sinkovic,
\emph{A graph for which the inertia bound is not tight},
J. Algebraic Combin. \textbf{47} (2018), 39--50.
\href{https://doi.org/10.1007/s10801-017-0768-0}
{doi:10.1007/s10801-017-0768-0}.

\bibitem{SS05}
N.~Srebro and A.~Shraibman,
\emph{Rank, trace-norm and max-norm},
in \emph{Learning Theory---COLT 2005},
Lecture Notes in Computer Science \textbf{3559},
Springer, 2005, pp.~545--560.
\href{https://doi.org/10.1007/11503415_37}
{doi:10.1007/11503415\_37}.

\bibitem{TZE}
Q.~Tang, S.~Zhang, and C.~Elphick,
\emph{Inertia, Independence and Expanders},
Bull. Lond. Math. Soc. \textbf{57} (2025), 4076--4095.
\href{https://doi.org/10.1112/blms.70185}{doi:10.1112/blms.70185}.
\url{https://arxiv.org/abs/2505.07305}.
\end{thebibliography}
\end{document}